\documentclass[11pt,
]{amsart}
\usepackage{xcolor}
\definecolor{winered}{rgb}{0.6,0,0}
\definecolor{lessblue}{rgb}{0,0,0.7}
\usepackage[pdftex,colorlinks=true,linkcolor=winered,citecolor=lessblue,breaklinks=true,bookmarksopen=true]{hyperref}
\usepackage{
amssymb, graphicx
}
\usepackage{mathrsfs} 
\usepackage{amsthm}
\usepackage{graphicx,color}
\newtheorem{theorem}{Theorem}
\newtheorem{lemma}{Lemma}
\newtheorem{proposition}{Proposition}
\newtheorem{assumption}{Assumption}

\theoremstyle{remark}

\date{\today}

\title[An IBVP for a semilinear wave equation]{An inverse boundary value problem for a semilinear wave equation on Lorentzian manifolds}

\author[P. Hintz]{Peter Hintz}
  \address{Department of Mathematics, Massachusetts Institute of Technology, Cambridge, MA 02139,USA
  (\tt{phintz@mit.edu})
  }
  \author[G. Uhlmann]{Gunther Uhlmann}
\address{Department of Mathematics, University of Washington, Seattle, WA 98195, USA; Institute for Advanced Study, 
The Hong Kong University of Science and Technology, Kowloon, Hong Kong, China (\tt{gunther@math.washington.edu})
}
  
  \author[J. Zhai]{Jian Zhai}
\address{Institute for Advanced Study,
  The Hong Kong University of Science and Technology, Kowloon, Hong Kong, China
  (\tt{iasjzhai@ust.hk}).}
 
\begin{document}
\begin{abstract}
We consider an inverse boundary value problem for a semilinear wave equation on a time-dependent Lorentzian manifold with time-like boundary. The time-dependent coefficients of the nonlinear terms can be recovered in the interior from the knowledge of the Neumann-to-Dirichlet map. Either distorted plane waves or Gaussian beams can be used to derive uniqueness.
\end{abstract}
\maketitle

\section{Introduction}
Let $(M,g)$ be a $(1+3)$-dimensional Lorentzian manifold with boundary $\partial M$, where the metric $g$ is of signature $(-,+,+,+)$. We assume that $M=\mathbb{R}\times N$ where $N$ is a manifold with boundary $\partial N$, and write the metric $g$ as
\begin{equation}
\label{EqIntroMetric}
g =-\beta(t,x')\mathrm{d}t^2+\kappa(t,x'),
\end{equation}
where $x=(t,x')=(x^0,x^1,x^2,x^3)$ are local coordinates on $M$; here, $\beta:\mathbb{R}\times N\rightarrow(0,\infty)$ is a smooth function and $\kappa(t,\cdot)$ is a Riemannian metric on $N$ depending smoothly on $t\in \mathbb{R}$. The boundary $\partial M=\mathbb{R}\times\partial N$ of $M$ is then timelike. Let $\nu$ denote the unit outer normal vector field to $\partial M$. Assume that $\partial M$ is null-convex, which means that $\mathrm{II}(V,V)=g(\nabla_V\nu,V)\geq 0$ for all null vectors $V\in T(\partial M)$; see \cite{hintz2017reconstruction} for a discussion of this condition. We consider the semilinear wave equation on $M$
\begin{alignat}{2}
\square_g u(x)+H(x,u(x))&=0,&\quad&\text{on }M,\nonumber\\
\label{maineq}
\partial_\nu u(x)&=f(x),&\quad& \text{on }\partial M,\\
u(t,x')&=0, &\quad& t<0,\nonumber
\end{alignat}
where $\square_g=|\det g|^{-1/2}\partial_j(\sqrt{|\det g|}g^{j k}\partial_k)$ is the wave operator (d'Alembertian) on $(M,g)$. We assume that $H(x,z)$ is smooth in $z$ near $0$ with Taylor expansion
\[
H(x,z) \sim \sum_{k=2}^\infty h_k(x)z^k,\quad h_k\in\mathcal C^\infty(M).
\]
As Neumann data, we take $f$ which are small in $\mathcal C^{m+1}$ for fixed large $m$. The Neumann-to-Dirichlet (ND) map $\Lambda$ is defined as
\[
\Lambda f = u\vert_{\partial M},
\]
where $u$ is the solution of \eqref{maineq}. We will investigate the inverse problem of determining $h_j(x)$, $j=2,3,\dots$, from $\Lambda$.

We remark that for the \emph{linear} equation $\square_g u+Vu=0$, the problem of recovering $V$ from the ND map is still open in general. Stefanov and Yang \cite{stefanov2018inverse} proved that the light ray transform of $V$ can be recovered from boundary measurements; however, the invertibility of the light ray transform is still unknown on general Lorentzian manifolds. We refer to \cite{lassas2019light,feizmohammadi2019light,vasy2019light} for an overview and recent results on the light ray transform.
 
In \cite{kurylev2018inverse}, the nonlinearity was exploited to solve inverse problems for a nonlinear equation where the corresponding inverse problem is still open for linear equations. The starting point of the approach is the higher order linearization, which we shall briefly introduce here. We take boundary Neumann data of the form $f=\sum_{i=1}^N\epsilon_if_i$, where $\epsilon_i$, $i=1,\ldots, N$ are small parameters. Since $\Lambda$ is a nonlinear map, $\Lambda(\sum_{i=1}^N\epsilon_if_i)$ contains more information than $\{\Lambda(f_i)\}_{i=1,\ldots, N}$: indeed, useful information can be extracted from
\[
\frac{\partial^N}{\partial\epsilon_1\cdots\partial\epsilon_N}\Big\vert_{\epsilon_1=\cdots=\epsilon_N=0}\Lambda\biggl(\sum_{i=1}^N\epsilon_if_i\biggr).
\]
This higher order linearization technique has been extensively used in the literature \cite{sun1997inverse,kang2002identification,kurylev2018inverse,lassas2018inverse,kurylev2014inverse,lassas2017determination,de2019nonlinear,wang2019inverse,de2018nonlinear,uhlmann2018determination,chen2019detection,assylbekov2017direct,carstea2019reconstruction,lassas2019inverse,lassas2019partial,feizmohammadi2019inverse,krupchyk2019partial,krupchyk2020remark,lai2020reconstruction}

The recovery of nonlinear terms from source-to-solution map was considered in \cite{lassas2018inverse}, where the authors use the nonlinear interactions of distorted plane waves. The approach originated from \cite{kurylev2018inverse}, and has been successfully used to study inverse problems for nonlinear hyperbolic equations \cite{lassas2018inverse,kurylev2014inverse,lassas2017determination,de2019nonlinear,wang2019inverse,de2018nonlinear,uhlmann2018determination,chen2019detection}. For some similar problems, Gaussian beams are used instead of distorted plane waves  \cite{kurylev2013determination,feizmohammadi2019recovery,uhlmann2019inverse}. The two approaches are actually closely related; both enable a pointwise recovery of the coefficients in the interior.

In this article, we will study the above inverse boundary value problem using both distorted plane waves and Gaussian beams. The two approaches will be discussed and compared in the last section.

To state our main result, recall that a smooth curve $\mu:(a,b)\rightarrow M$ is \emph{causal} if $g(\dot{\mu}(s),\dot{\mu}(s))\leq 0$ and $\dot{\mu}(s)\neq 0$ for all $s\in(a,b)$. Given $p,q\in M$, we write $p\leq q$ if $p= q$ or $p$ can be joined to $q$ by a future directed causal curve. We say $p<q$ if $p\leq q$ and $p\neq q$. We denote the causal future of $p\in M$ by $J^+(p)=\{q\in M: p\leq q\}$ and the causal past of $q\in M$ by $J^-(q)=\{p\in M: p\leq q\}$. We shall restrict the ND map to $(0,T)\times\partial N$, and correspondingly work in
\[
\mathbb{U}=\bigcup_{p,q\in (0,T)\times\partial N}J^+(p)\cap J^-(q).
\]

\begin{theorem}\label{maintheorem}
Consider the semilinear wave equations
\[
\square_gu(x)+H^{(j)}(x,u(x))=0,\qquad j=1,2.
\]
Assume $H^{(j)}(x,z)$ are smooth in $z$ near $0$ and have a Taylor expansion\footnote{The notation means that $h_k^{(j)}(x)=\frac{1}{k!}\frac{\partial^k}{\partial z^k}H^{(j)}(x,0)$.}
\[
H^{(j)}(x,z)\sim\sum_{k=2}^\infty h_k^{(j)}(x)z^k,\qquad h_k^{(j)}\in\mathcal C^\infty(\overline{\mathbb{U}}). 
\]
Assume that null geodesics in $\mathbb{U}$ do not have cut points. If the Neumann-to-Dirichlet maps $\Lambda^{(j)}$ acting on $\mathcal C^6([0,T]\times\partial N)$ are equal, $\Lambda^{(1)}=\Lambda^{(2)}$, then
\[
  h^{(1)}_k(x)=h_k^{(2)}(x),\qquad x\in\mathbb{U},\ k\geq 2.
\]
\end{theorem}

The strategy of the proof is to send in distorted plane waves (or Gaussian beams) from \emph{outside} the manifold $M$ (within a small extension $\widetilde M$) and analyze contributions to the ND map from nonlinear interactions in the interior of $M$ as well as from subsequent reflections at the boundary $\partial M$ of $M$.

The rest of this paper is organized as follows. In Section \ref{wellposedness}, we establish the well-posedness of the initial boundary value problem \eqref{maineq} for small boundary data. In Section \ref{distorted}, we use the nonlinear interaction of distorted plane waves to prove the main theorem. In Section \ref{gaussianbeam}, we give another proof of the main theorem using Gaussian beam solutions, assuming $h_2$ is already known. Finally, the two approaches will be compared and discussed in Section \ref{discussion}.

\section{Well-posedness for small boundary data}\label{wellposedness}
We establish well-posedness of the initial boundary value problem \eqref{maineq} in this section with small boundary value $f$.

Fix $m\geq 5$.  We assume $f\in \mathcal C^{m+1}([0,T]\times\partial N)$ and $\|f\|_{\mathcal C^{m+1}([0,T]\times\partial N)}\leq \epsilon_0$
 for a small number $\epsilon_0>0$. Assume also that $f$ satisfies the compatibility condition $\frac{\partial^\ell f}{\partial t^\ell}=0$ at $\{t=0\}$ for any $\ell=0,1,\dots,m-1$. We can find a function $h\in \mathcal C^{m+1}([0,T]\times N)$ such that $\partial_\nu h\vert_{[0,T]\times\partial N}=f$ and
\[
\|h\|_{\mathcal C^{m+1}([0,T]\times N)}\leq C\|f\|_{\mathcal C^{m+1}([0,T]\times\partial N)}.
\]
Let $\widetilde{u}=u-h$, where $u$ solves the initial boundary value problem \eqref{maineq}. Then $\widetilde{u}$ satisfies the equation
\[
\square_g\widetilde{u}=F(x,\widetilde{u},h) := -\square_g h-H(x,\widetilde{u}+h),
\]
supplemented with the boundary condition $\partial_\nu\widetilde{u}=0$ on $(0,T)\times\partial N$ and initial conditions $\widetilde{u}=\frac{\partial\widetilde{u}}{\partial t}=0$ at $\{0\}\times N$.
The above equation can be written in the form
\begin{alignat}{2}
&\square_g\widetilde{u}=F(x,\widetilde{u},h),&\quad&\text{ in }(0,T)\times\partial N,\nonumber\\
\label{eq_hom}
&\partial_\nu\widetilde{u}=0,&\quad&\text{ on }(0,T)\times\partial N,\\
&\widetilde{u}=\frac{\partial\widetilde{u}}{\partial t}=0,&\quad&\text{ on }t=0. \nonumber
\end{alignat}
This equation is of the form \cite[equation (5.12)]{dafermos1985energy}. For $R>0$, define $Z(R, T)$ as the set of all functions $w$ satisfying
\[
w\in \bigcap_{k=0}^m W^{k,\infty}([0,T];\, H^{m-k}(N)),\qquad
\|w\|_Z^2 := \sup_{t\in[0,T]}\sum_{k=0}^m\|\partial^k_tw(t)\|^2_{H^{m-k}}\leq R^2.
\]
We can write $F(x,\widetilde{u},h)=\mathscr{F}(x,h)+G(x,\widetilde{u},h)\widetilde{u}$ where $\mathscr{F}=-\square_gh-H(x,h)$ and
\[
G(x,\widetilde{u},h)=-\int_0^1\partial_zH(x,h+\tau\widetilde{u})\mathrm{d}\tau.
\]
We can write $\mathscr{F}(x,h)=\mathscr{F}(t,y,h)$ using the notation $x=(t,y)$. Since $H(x,z)$ is smooth in $z$, we have
\[
  \sup_{t\in[0,T]}\sum_{k=0}^{m-1}\|\partial^{k}_t\mathscr{F}(t,\cdot,h)\|_{H^{m-k-1}}\leq C\sup_{t\in[0,T]}\sum_{k=0}^{m-1}\|\partial^k_t\mathscr{F}(t,\cdot,h)\|_{\mathcal C^{m-k-1}} \leq C'\epsilon_0.
\]
Moreover, $\partial_z H(x,z)$ vanishes linearly in $z$, hence we have
\[
  G(x,\widetilde{u},h)\in\bigcap_{k=0}^m  W^{k,\infty}([0,T];\, H^{m-k}(N)), \quad
  \|G(x,\widetilde{u},h)\|_Z \leq C(\|h\|_Z+\|\widetilde{u}\|_Z) \leq C'(\epsilon_0+\|\widetilde{u}\|_Z)
\]
for $\widetilde{u}\in Z(\rho_0,T)$ with $\rho_0$ small enough.

Given $\widetilde{w}\in Z(\rho_0,T)$, consider first the \emph{linear} initial boundary value problem
\begin{alignat}{2}
&\square_g\widetilde{u}-G(x,\widetilde{w},h)\widetilde{w}=\mathscr{F}(x,h),&\quad& t\in(0,T), \nonumber\\
\label{linearized_IVP}
&\partial_\nu\widetilde{u} = 0,&\quad& t\in(0,T), \\
&\widetilde{u}(0)=\frac{\partial \widetilde{u}}{\partial t}(0)=0.\nonumber
\end{alignat}
By \cite[Theorem 3.1]{dafermos1985energy}, there exists a unique solution $\widetilde{u}\in \bigcap_{k=0}^m \mathcal C^k([0,T];\, H^{m-k}(N))$
to \eqref{linearized_IVP}, and it satisfies the estimate
\[
  \|\widetilde{u}\|_Z \leq C(\epsilon_0 + \epsilon_0\|\widetilde{w}\|_Z+\|\widetilde{w}\|_Z^2) e^{K T},
\]
where $C,K$ are positive constants depending on the coefficients of the equation. Denote $\mathscr{T}$ to be the map which maps $\widetilde{w}\in Z(\rho_0,T)$ to the solution $\widetilde{u}$ of \eqref{linearized_IVP}. Notice that we can take $\rho_0$ small enough and $\epsilon_0=\frac{e^{-KT}}{2C}\rho_0$ such that
\[
  C(\epsilon_0+\epsilon_0\rho_0+\rho_0^2)e^{K T} < \rho_0.
\]
Then $\mathscr{T}$ maps $Z(\rho_0,T)$ to itself.

Now assume $\widetilde{u}_j$, $j=1,2$, solve the equation
\[
\begin{split}
&\square_g\widetilde{u}_j-G(x,\widetilde{w}_j,h)\widetilde{w}_j=\mathscr{F}(x,h),\quad t\in(0,T)\\
&\widetilde{u}_j(0)=\frac{\partial \widetilde{u}_j}{\partial t}(0)=0.
\end{split}
\]
We have $\widetilde{u}_j=\mathscr{T}\widetilde{w}_j$, $j=1,2$ and
\[
\square_g(\widetilde{u}_1-\widetilde{u}_2)=-\left(\int_0^1\partial_zH(x,h+\widetilde{w}_2+\tau(\widetilde{w}_1-\widetilde{w}_2))\mathrm{d}\tau\right)(\widetilde{w}_1-\widetilde{w}_2).
\]
Then
\[
  \|\mathscr{T}\widetilde{w}_1-\mathscr{T}\widetilde{w}_2\|_Z=\|\widetilde{u}_1-\widetilde{u}_2\|_Z\leq C (\epsilon_0+\rho_0) e^{K T}\|\widetilde{w}_1-\widetilde{w}_2\|_Z.
\]
Choosing $\rho_0$ small enough such that $C (\epsilon_0+\rho_0)e^{K T}<1$, the map $\mathscr{T}$ is a contraction. Consequently, the equation \eqref{eq_hom} has a unique solution $\widetilde{u}$ in $Z(\rho_0,T)$. Using
 \cite[Theorem 3.1]{dafermos1985energy} again, we have
 \[
 \widetilde{u}\in\bigcap_{k=0}^m\mathcal{C}^k([0,T];\, H^{m-k}(N)).
 \]
In summary, we have shown:

\begin{theorem}
Let $T>0$ be fixed. Assume that $f\in \mathcal C^{m+1}([0,T)\times\partial N)$, $m\geq 5$, and $\frac{\partial^\ell f}{\partial t^\ell}=0$ at $\{t=0\}$ for any $\ell=0,1\cdots, m-1$ at $t=0$. Then there exists $\epsilon_0>0$ such that for $\|f\|_{\mathcal C^{m}}\leq \epsilon_0$, there exists a unique solution
\[
u\in \bigcap_{k=0}^m\mathcal{C}^k([0,T];\, H^{m-k}(N))
\]
of equation \eqref{maineq}. It satisfies the estimate
\[
   \sup_{t\in[0,T]}\|\partial^{m-k}_tu(t)\|_{H^{m-k}(N)}\leq C\|f\|_{\mathcal C^{m+1}([0,T]\times\partial N)}, 
\]
where $C>0$ is independent of $f$.

\end{theorem}

If $f=\epsilon f_1$ where $\epsilon>0$ is small, then for any $N=1,2,\dots$, we can write (cf. \cite[Appendix III]{choquet2008general} and the discussion in \cite[Section 3.1]{kurylev2018inverse})
\begin{equation}
\label{EqAnsatz}
u=\sum_{j=1}^N\epsilon^jw_j+R_N,
\end{equation}
where $w_j\in \bigcap_{k=0}^m\mathcal{C}^k([0,T];\, H^{m-k}(N))
$ for $j=1,\cdots, N$, $R_N\in \bigcap_{k=0}^m\mathcal{C}^k([0,T];\, H^{m-k}(N))$ and
\[
\sup_{t\in[0,T]}\|\partial^{m-k}_tR_N(t)\|_{H^{m-k}(N)}\leq C_N\epsilon^{N+1},
\]
where $C_N>0$ is a constant depending on $N$. Indeed, this follows by plugging~\eqref{EqAnsatz} as an ansatz into equation~\eqref{maineq}, solving inductively for the coefficients $w_j$ (which only involves the solution of linear wave equations), and solving a nonlinear equation for $R_N$ with forcing term of size $\epsilon^{N+1}$. Hence one can denote
\begin{equation}\label{def_deriv}
w_N=\frac{\partial^N}{\partial\epsilon^N}u\vert_{\epsilon=0}.
\end{equation}
The proof presented later will heavily depend on the above asymptotic expansion.

\section{Recovery using distorted plane waves}\label{distorted}
In this section we will show how to recover $h_k$, $k=2,3,\ldots$ by using the nonlinear interaction of distorted plane waves.
First we extend the metric $g$ on $M$ smoothly to a metric $\widetilde{g}$ on a larger manifold $\widetilde{M}=\mathbb{R}_t\times\widetilde{N}$ such that
\begin{enumerate}
\item $N$ is contained in the interior of $\widetilde{N}$, and thus $M$ is contained in the interior of $\widetilde{M}$;
\item $\widetilde{N}$ is closed, i.e.\ compact without boundary,
\item $\widetilde{g}$ is a warped product metric, $\widetilde{g}=-\widetilde{\beta}(t,x')\mathrm{d}t^2+\widetilde{\kappa}(t,x')$, with $\widetilde{\beta}=\beta$ and $\widetilde{\kappa}=\kappa$ on $M$ in the notation of~\eqref{EqIntroMetric}.
\end{enumerate}
We can for example take $\widetilde{N}$ to be the double of $N$, and define $\widetilde{\beta}$ to be an arbitrary but smooth and positive extension of $\beta$ to $\widetilde{M}$, and similarly $\widetilde{\kappa}$ to be an arbitrary but smooth positive section of $S^2 T^*\widetilde{N}$ over $\mathbb{R}_t\times\widetilde{N}$ extending $\kappa$. The advantage of this construction is that $\widetilde{M}$ is globally hyperbolic, which will occasionally be useful.

\subsection{Notations and preliminaries}\label{notation}

For $p\in\widetilde M$, denote the set of light-like vectors at $p$ by
\[
L_p\widetilde M=\{\zeta\in T_p\widetilde M\setminus\{0\}:\, g(\zeta,\zeta)=0\}.
\]
The set of light-like covectors at $p$ is denoted by $L^{*}_p\widetilde M$. The sets of future and past light-like vectors (covectors) are denoted by $L_p^+\widetilde M$ and $L_p^-\widetilde M$ ($L^{*,+}_p\widetilde M$ and $L^{*,+}_p\widetilde M$).
Define the future directed light-cone emanating from $p$ by
\[
\mathcal{L}^+(p)=\{\gamma_{p,\zeta}(t)\in \widetilde M:\zeta\in L^+_p\widetilde M, t\geq 0\}\subset \widetilde M.
\]

Distorted plane waves have singularities conormal to a submanifold of $\widetilde M$ and can be viewed as Lagrangian distributions. We review them briefly, closely following the notation used in \cite{lassas2018inverse}. Recall that $T^*\widetilde M$ is a symplectic manifold with canonical 2-form, given in local coordinates by $\omega=\sum_{j=1}^4\mathrm{d}\xi_j\wedge\mathrm{d}x^j$. A submanifold $\Lambda\subset T^*\widetilde M$ is called Lagrangian if $n:=\mathrm{dim}\, \Lambda=4$ and $\omega$ vanishes on $\Lambda$. For $K$ a smooth submanifold of $\widetilde M$, its conormal bundle
\[
N^*K=\{(x,\zeta)\in T^*\widetilde M: x\in K,\, \langle \zeta,\theta\rangle=0, \,\theta\in T_xK\}
\]
is a Lagrangian submanifold of $T^*\widetilde M$. Let $\Lambda$ be a smooth conic Lagrangian submanifold of $T^*\widetilde M\setminus 0$. We denote by $\mathcal{I}^\mu(\Lambda)$ the space of Lagrangian distributions of order $\mu$ associated with $\Lambda$. If $\Lambda=N^*K$ for some submanifold $K\subset\widetilde{M}$, then $\mathcal{I}^\mu(K):=\mathcal{I}^\mu(N^*K)$ denotes the space of conormal distributions to $K$. For $u\in\mathcal{I}^\mu(\Lambda)$, one can define the principal symbol $\sigma^{(p)}(u)=\sigma^{(p)}_{\Lambda}(u)$ of $u$ with
\[
\sigma^{(p)}(u)\in S^{\mu+\frac{n}{4}}(\Lambda,\Omega^{1/2}\otimes L)/S^{\mu+\frac{n}{4}-1}(\Lambda,\Omega^{1/2}\otimes L),
\]
where $\Omega^{1/2}$ is the half-density on $\widetilde{M}$ and $L$ is the Maslov--Keller line bundle of $\Lambda$. We refer to \cite[Chapter 4]{duistermaat1996fourier} for the precise definition and more discussions.

For waves described by nonlinear wave equations, the distorted plane waves, characterized by Lagrangian distributions, can have nonlinear interactions and generate new propagating singularities. Such new singularities can be characterized by paired Lagrangian distributions, which will be reviewed below. The detailed analysis of the singularities and principal symbols of the waves generated by nonlinear interactions is the key to the study of various inverse problems for nonlinear wave equations \cite{kurylev2018inverse, lassas2018inverse,kurylev2014inverse,lassas2017determination,de2019nonlinear,wang2019inverse,de2018nonlinear,uhlmann2018determination,chen2019detection}. Let $\Lambda_0,\Lambda_1\subset T^*\widetilde M\setminus 0$ be two Lagrangian submanifolds intersecting cleanly, i.e.,
\[
T_p\Lambda_0\cap T_p\Lambda_1=T_p(\Lambda_0\cap\Lambda_1)\quad \forall\,p\in \Lambda_0\cap \Lambda_1.
\]
We denote the space of paired Lagrangian distributions associated with $(\Lambda_0,\Lambda_1)$ by $\mathcal{I}^{p,l}(\Lambda_0,\Lambda_1)$.  We mention here that if $u\in\mathcal{I}^{p,l}(\Lambda_0,\Lambda_1)$, then microlocally away from $\Lambda_0\cap\Lambda_1$, we have $u\in \mathcal{I}^{p+l}(\Lambda_0\setminus\Lambda_1)$ and $u\in \mathcal{I}^{p}(\Lambda_1\setminus\Lambda_0)$ with well defined principal symbols $\sigma^{(p)}_{\Lambda_0}(u)$ and $\sigma^{(p)}_{\Lambda_1}(u)$. For more details, we refer to \cite{melrose1979lagrangian,guillemin1981oscillatory}.

Fix a Riemannian metric $g^+$ on $\widetilde M$. Given $x_0\in\widetilde{M}\setminus M$, $\zeta_0\in L^+_{x_0}\widetilde{M}$, and $s_0>0$, put
\begin{align*}
\mathcal{W}_{x_0,\zeta_0,s_0}&=\{\eta\in L_{x_0}^+\widetilde{M}:\,\|\eta-\zeta_0\|_{g^+}<s_0,\,\|\eta\|_{g^+}=\|\zeta_0\|_{g^+}\}, \\
K(x_0,\zeta_0,s_0)&=\{\gamma_{x_0,\eta}(s)\in \widetilde{M}:\,\eta\in \mathcal{W}_{x_0,\zeta_0,s_0},s\in(0,\infty)\}, \\
\Lambda(x_0,\zeta_0,s_0)&=\{(\gamma_{x_0,\eta}(s),r\dot{\gamma}_{x_0,\eta}(s)^\flat)\in T^*\widetilde{M};\,\eta\in \mathcal{W}_{x_0,\zeta_0,s_0},s\in (0,\infty),r>0\}.
\end{align*}
Notice that $K(x_0,\zeta_0,s_0)$ is a subset of codimension $1$ of the light cone $\mathcal{L}^+(x_0)$, and
\[
N^*K(x_0,\zeta_0,s_0)=\Lambda(x_0,\zeta_0,s_0).
\]
By \cite[Lemma 3.1]{kurylev2018inverse}, one can construct distributions $u_0\in\mathcal{I}^{\mu}(\widetilde M\setminus \{x_0\},\Lambda(x_0,\zeta_0,s_0))$ which on $M$ satisfy $\square_g u_0\in\mathcal C^\infty(M)$, and whose principal symbol is nonzero on $(\gamma_{x_0,\zeta_0}(s),\dot\gamma_{x_0,\zeta_0}(s)^\flat)$. Thus, $u_0$ is a nontrivial distorted plane wave propagating on the surface $K(x_0,\zeta_0,s_0)$.

We consider four distorted plane waves
\[
u_j\in\mathcal{I}^{\mu}(\widetilde M,\Lambda(x_j,\xi_j,s_0)),\quad j=1,2,3,4,
\]
which are approximate solutions of the linearized wave equation in $M$, that is, $\square_g u_j\in\mathcal{C}^\infty(M)$. Let
\begin{equation}\label{KLambda}
K_j=K(x_j,\xi_j,s_0),\quad\quad \Lambda_j=\Lambda(x_j,\xi_j,s_0)=N^*K_j.
\end{equation}
As in \cite{lassas2018inverse}, we make the following assumptions.
\begin{assumption}\label{assumption1}
Assume that 
\begin{enumerate}
\item $K_i,\,K_j$, $i\neq j$, intersect at a codimension $2$ submanifold $K_{ij}\subset\widetilde{M}$;
\item $K_i,\, K_j,\, K_k$, $i,j,k$ distinct, intersect at a codimension $3$ submanifold $K_{ijk}\subset\widetilde{M}$;
\item $K_1,\,K_2,\,K_3,\,K_4$ intersect at a point $q_0\in M$.
\end{enumerate}
Assume further that for any two disjoint subsets $I,J\subset\{1,2,3,4\}$, the intersection of $\cap_{i\in I}K_i$ and $\cap_{j\in J}K_j$ is transversal if not empty.
\end{assumption}

We use the notations
\[
 \Lambda_{ij}=N^*K_{ij},\quad \Lambda_{ijk}=N^*K_{ijk},\qquad
 \Lambda_{q_0}=T^*_{q_0}M\setminus 0;
\]
which are all Lagrangian submanifolds in $T^*\widetilde{M}$.
For any $\Gamma\subset T^*\widetilde{M}$, we denote by $\Gamma^g$ the flow-out of $\Gamma\cap L^{*,+}\widetilde M$ under the null-geodesic flow of $g$ lifted to $T^*\widetilde{M}$. To define this precisely, denote by $H_G\in\mathcal C^\infty(T^*\widetilde{M};T T^*\widetilde{M})$ the Hamilton vector field of the dual metric function $G\colon T^*\widetilde{M}\ni\zeta\mapsto g^{-1}(\zeta,\zeta)$; in particular $L^*\widetilde{M}=G^{-1}(0)$. We then put
\begin{equation}
\label{EqFlowout}
  \Gamma^g:=\{ \exp(s H_G)\zeta \colon 0\leq s<s_+(\zeta),\ \zeta\in\Gamma\cap L^{*,+}\widetilde{M}\},
\end{equation}
where $s\mapsto\exp(s H_G)\zeta$ is the integral curve of $H_G$ with initial condition $\zeta$, and $s_+(\zeta)\in(0,\infty)\cup\{+\infty\}$ is the supremum of the maximal interval of existence of this integral curve.

We assume $x_j\in (0,T)\times\widetilde{N}$; we can take $s_0$ small enough so that $u_j$ is smooth near $t=0$. Denote $f_i=\partial_\nu u_i\vert_{\partial M}$; then the solutions $v_i$ of the linear equations
\begin{alignat*}{2}
\square_g v_i(x)&=0,&\quad&\text{on }M,\\
\partial_\nu v_i(x)&=f_i(x),&\quad& \text{on }\partial M,\\
v_i(t,x')&=0, &\quad& t<0,
\end{alignat*}
are equal to $u_i$ modulo $\mathcal C^\infty(M)$. We will use nonlinear interactions of three or four distorted plane waves for our study. For $N=3$ or $4$, consider then
\begin{equation}\label{boundarysource}
f=\sum_{i=1}^N\epsilon_i f_i,
\end{equation}
and denote $v=\sum_{i=1}^N \epsilon_iv_i$. We write $w=Q_g(F)$ if $w$ solves the linear wave equation
\begin{equation}\label{def_Q}
\begin{split}
\square_g w(x)&=F,\quad\text{on }M,\\
\partial_\nu w(x)&=0,\quad \text{on }\partial M,\\
w&=0, \quad t<0.
\end{split}
\end{equation}
The solution $u$ to \eqref{maineq} is then given by the asymptotic expansion \cite[(2.9)]{lassas2018inverse}
\begin{equation}\label{asymtotic_u}
\begin{split}
u=v&-Q_g(h_2v^2)+2Q_g(h_2vQ_g(h_2v^2)-4Q_g(h_2vQ_g(h_2vQ_g(h_2v^2)))\\
&-Q_g(h_2Q_g(h_2v^2)Q_g(h_2v^2))+2Q_g(h_2vQ_g(h_3v^3))-Q_g(h_3v^3)+3Q_g(h_3v^2Q_g(h_2v^2))\\
&-Q_g(h_4v^4)+\text{higher order terms in } \epsilon_1,\ldots,\epsilon_N.
\end{split}
\end{equation}
We will use the singularities from the terms in \eqref{asymtotic_u} to recover the coefficients of \eqref{maineq}. Notice that those terms involve nonlinear interactions of distorted plane waves $v_j$, $j=1,\dots,N$, and thus new singularities can be created. Recovery of a Lorentzian metric from the source-to-solution map using those newly generated singularities was first carried out in \cite{kurylev2018inverse}. For recovery of the coefficients of nonlinear terms, we refer to \cite{lassas2018inverse,de2019nonlinear}.

\subsection{Nonlinear interactions of three waves and recovery of \texorpdfstring{$(h_2)^2$ and $h_3$}{h3 and the square of h2}}\label{threewaves}
First, we will first use three distorted plane waves, i.e.\ taking $N=3$ in \eqref{boundarysource} and using Neumann data
\[
f=\sum_{i=1}^3\epsilon_i f_i
\]
with $\epsilon_i>0$, $i=1,2,3$, small parameters. We will construct suitable sources $f_i$, $i=1,2,3$, and denote by $v_i$ the corresponding distorted plane wave.

For any $p\in M$ and $\xi\in L^{*,+}_pM$ define $\gamma(s)=\gamma_{p,\xi}(s)$ to be the geodesic such that $\gamma(0)=p$ and $\dot{\gamma}(0)=\xi^\sharp$. Define
\[
s^+(p,\xi)=\inf\{s>0:\gamma(s)\in \partial M\},\quad s^-(p,\xi)=\sup\{s<0:\gamma(s)\in \partial M\}.
\]
Fix a point $q_0\in \mathbb{U}$. There exist $\xi^{(0)},\xi^{(1)}\in L^{*,+}_{q_0}M$ such that
\begin{equation}
\label{Eqxminus}
x^-=\gamma_{q_0,\xi^{(1)}}(s^-(q_0,\xi^{(1)}))\in (0,T)\times\partial N,\quad x_0=\gamma_{q_0,\xi^{(0)}}(s^+(q_0,\xi^{(0)}))\in (0,T)\times\partial N.
\end{equation}
Indeed, by definition of $\mathbb{U}$, there exists a point $(t_0,y_0)\in(0,T)\times\partial N$ and a future causal curve lying inside $M$ which joins $(t,y)$ and $q_0$. Since the $t$-coordinate of $q_0$ is less than $T$, the set of $t\in[t_0,T)$ for which there exists a future causal curve inside the larger manifold $\widetilde{M}$ joining $(t,y_0)$ and $q_0$ has a least upper bound $\bar t<T$. Standard compactness arguments on the globally hyperbolic manifold $\widetilde{M}$ imply that there exists a future causal curve $\gamma$ from $(\bar t,y_0)$ to $q_0$, which by short-cut arguments must be a positive reparameterization of a null-geodesic without cut points \cite[\S10]{ONeillSemi}. Upon normalizing $\gamma$ so that $\gamma(0)=(\bar t,y_0)$ and $\gamma(1)=q_0$, the backwards null-geodesic $\mu\colon[0,s_0]\to\widetilde{M}$ with initial data $(q_0,-\dot\gamma(1))$ coincides with $\gamma$ until it reaches $\mu(s_0)=(\bar t,y_0)$. Note that $t\circ\mu\colon[0,s_0]\to(0,T)$ is monotonically decreasing; thus, for the smallest $s\in(0,s_0]$ so that $x^-=\mu(s)\in(0,T)\times\partial N$, we necessarily have $(t\circ\mu)(s)\in[(t\circ\mu)(s_0),(t\circ\mu)(0)]\subset(0,T)$. This shows that $x^-\in(0,T)\times\partial N$ is of the form~\eqref{Eqxminus} with $\xi^{(1)}=-\dot\mu(s)^\flat$, and in particular proves the existence of $\xi^{(1)}$. The argument for $\xi^{(0)}$ is analogous.

Put $\gamma^{(j)}=\gamma_{q_0,\xi^{(j)}}$, $j=0,1$ and denote $x_1=\gamma^{(1)}(s^-(q_0,\xi^{(1)})-\epsilon)$ for $\epsilon>0$ small; thus, $x_1\in\widetilde{M}\setminus M$ lies just barely outside of $M$. 

Choose local coordinates so that $g$ coincides with the Minkowski metric at $q_0$. Using further linear changes of coordinates which leave the Minkowski metric unchanged (that is, rotations in the spatial variables, Lorentz boosts), and upon scaling $\xi^{(0)},\xi^{(1)}$ by a positive scalar, one can assume without loss of generality (cf. \cite[Lemma 1]{chen2019detection}) that
\[
\xi^{(0)}=(-1,-\sqrt{1-r_0^2},r_0,0),\quad\quad \xi^{(1)}=(-1,1,0,0),
\]
for some $r_0\in [-1,1]$. Take a small parameter $\varsigma>0$ and introduce two perturbations of $\xi^{(1)}$
\[
\xi^{(2)}=(-1,\sqrt{1-\varsigma^2},\varsigma,0),\quad\quad \xi^{(3)}=(-1,\sqrt{1-\varsigma^2},-\varsigma,0).
\]
Notice $\xi^{(2)},\xi^{(3)}\in L^{*,+}_pM$. One can then write $\xi^{(0)}$ as a linear combination of $\xi^{(1)},\xi^{(2)},\xi^{(3)}$,
\[
\xi^{(0)}=\alpha_1\xi^{(1)}+\alpha_2\xi^{(2)}+\alpha_3\xi^{(3)},
\]
with
\[
\alpha_1=\frac{-\sqrt{1-\varsigma^2}-\sqrt{1-r_0^2}}{1-\sqrt{1-\varsigma^2}},\quad \alpha_2=\frac{1+\sqrt{1-r_0^2}}{2(1-\sqrt{1-\varsigma^2})}+\frac{r_0}{2\varsigma},\quad \alpha_3=\frac{1+\sqrt{1-r_0^2}}{2(1-\sqrt{1-\varsigma^2})}-\frac{r_0}{2\varsigma}.
\]
Denote $b(r_0)=1+\sqrt{1-r_0^2}$. By direct calculation, and using the asymptotics $\sqrt{1-\varsigma^2}=1-\frac{1}{2}\varsigma^2+\mathcal{O}(\varsigma^4)$, we obtain
\begin{alignat*}{2}
|\alpha_1\xi^{(1)}+\alpha_2\xi^{(2)}|_g^2&=\ &2b(r_0)^2\varsigma^{-2}+\mathcal{O}(\varsigma^{-1}),\\
|\alpha_1\xi^{(1)}+\alpha_3\xi^{(3)}|_g^2&=\ &2b(r_0)^2\varsigma^{-2}+\mathcal{O}(\varsigma^{-1}),\\
|\alpha_2\xi^{(2)}+\alpha_3\xi^{(3)}|_g^2&=\ &-4b(r_0)^2\varsigma^{-2}+\mathcal{O}(\varsigma^{-1}).
\end{alignat*}
Therefore,
\begin{equation}
\label{nonvanishing30}
|\alpha_1\xi^{(1)}+\alpha_2\xi^{(2)}|_g^{-2}+|\alpha_1\xi^{(1)}+\alpha_3\xi^{(3)}|_g^{-2}+|\alpha_2\xi^{(2)}+\alpha_3\xi^{(3)}|_g^{-2}=\frac{3}{4b(r_0)^2}\varsigma^2+\mathcal{O}(\varsigma^3).
\end{equation}
By taking $\varsigma$ small enough, the quantity
\begin{equation}\label{nonvanishing3}
\sum_{\sigma\in\Sigma(3)}\left|\alpha_{\sigma(2)}\xi^{(\sigma(2))}+\alpha_{\sigma(3)}\xi^{(\sigma(3))}\right|_{g(q_0)}^{-2}
\end{equation}
is nonvanishing; here, $\Sigma(3)$ denotes the permutation group of $\{1,2,3\}$.

For $j=2,3$, let $\gamma^{(j)}=\gamma_{q_0,\xi^{(j)}}$, and denote
\[
x_j=\gamma^{(j)}(s^-(q_0,\xi^{(j)})-\epsilon),\quad j=2,3,
\]
for $\epsilon>0$ small. Here, if we took $\epsilon=0$, then we could choose $\varsigma$ small enough so that $x_j\in (0,T)\times\partial N$; fixing $\varsigma$ in this manner, we can then take $\epsilon>0$ small enough so that $x_j\in\widetilde{M}\setminus M$ and $t>0$ at $x_j$ still. Here we used the fact that null-geodesics are non-tangential, hence transversal, to $\partial M$ due to the null-convexity of $\partial M$. Now for $j=1,2,3$ denote
\[
\xi_j=\dot\gamma_{q_0,\xi^{(j)}}(s^-(q_0,\xi^{(j)})-\epsilon)^\flat\in L^{*,+}_{x_j}M.
\]
Use these $(x_j,\xi_j)$, $j=1,2,3$, in \eqref{KLambda} and denote associated distorted plane waves by
\[
u_j\in\mathcal{I}^{\mu}(\Lambda_j),\quad j=1,2,3.
\]
We note that $\xi^{(0)}\in N^*_pK_{123}$.

Let $u$ denote the solution of \eqref{maineq} with $f=\sum_{i=1}^3\epsilon_if_i$, and put 
\[
\mathcal{U}^{(3)}=\partial_{\epsilon_1}\partial_{\epsilon_2}\partial_{\epsilon_3}u\vert_{\epsilon_1=\epsilon_2=\epsilon_3=0},
\]
which can be defined analogous to \eqref{def_deriv}.
We can then decompose
\begin{equation}
\label{EqU3decomp}
\mathcal{U}^{(3)}=\mathcal{U}_0^{(3)}+\mathcal{U}_1^{(3)},\qquad
\mathcal{U}_0^{(3)}:=-6Q_g(h_3v_1v_2v_3),\quad \mathcal{U}_1^{(3)}:=2\sum_{\sigma\in\Sigma(3)}Q_g(h_2v_{\sigma(1)}Q(h_2v_{\sigma(2)}v_{\sigma(3)})).
\end{equation}
Recall that $Q_g$ is the solution operator associated with the equation \eqref{def_Q}. On globally hyperbolic manifolds (no boundary!), the wave operator $\square_g$ has a causal (retarded) inverse (cf. \cite[Theorem 3.3.1]{bar2007wave}). 
Denote by $\widetilde{Q}_g=\square_g^{-1}$ the causal inverse of $\square_g$ on $\widetilde{M}$. Then
\begin{equation}
\label{EqU3incDecomp}
\mathcal{U}^{(3),\mathrm{inc}}:=\mathcal{U}^{(3),\mathrm{inc}}_0+\mathcal{U}^{(3),\mathrm{inc}}_1=-6\widetilde{Q}_g(h_3v_1v_2v_3)+2\sum_{\sigma\in\Sigma(3)}\widetilde{Q}_g(h_2v_{\sigma(1)}\widetilde{Q}(h_2v_{\sigma(2)}v_{\sigma(3)}))
\end{equation}
is the incident wave before reflection on the boundary. We have (cf. \cite[Proposition 3.7]{lassas2018inverse} and the subsequent discussion):
\begin{proposition}\label{symbol_threewaves0}
Let $\Lambda^g_{123}$ be the flow-out of $\Lambda_{123}\cap L^{*,+}\widetilde M$, as defined in general in~\eqref{EqFlowout}. Then
\[
\mathcal{U}^{(3),\mathrm{inc}}\in \mathcal{I}^{3\mu+\frac{1}{2},-\frac{1}{2}}(\Lambda_{123},\Lambda^g_{123})
\]
away from $\cup_{i=1}^3\Lambda_i$. For any $q\in K_{123}$ and $\zeta\in N_{q}K_{123}$, there exists a unique decomposition $\zeta=\sum_{j=1}^3\zeta_j$ with $\zeta_j\in N_q^*K_j$ (cf. \cite{lassas2018inverse}). Assume that $(y,\eta)$ lies along the forward null-bicharacteristic of $\square_g$ starting at $(q,\zeta)$. The principal symbol of $\mathcal{U}^{(3),\mathrm{inc}}$ can be written as
\begin{equation}\label{U3}
\sigma^{(p)}(\mathcal{U}^{(3),\mathrm{inc}})(y,\eta)=\sigma^{(p)}(\mathcal{U}_0^{(3),\mathrm{inc}})(y,\eta)
=-6(2\pi)^{-2}\sigma^{(p)}(\widetilde{Q}_g)(y,\eta,q,\zeta)h_3(q)\prod_{j=1}^3\sigma^{(p)}(v_j)(q,\zeta_j).
\end{equation}
\end{proposition}

Note here that $K_{1 2 3}$ is a 1-dimensional spacelike submanifold since its conormal bundle at $p\in K_{1 2 3}$ is timelike, being the 1-codimensional vector space $N^*_p K_{1 2 3}=N^*_p K_1+N^*_p K_2+N^*_p K_3\subset T^*_p M$ which by assumption contains $3$ linearly independent null covectors. This implies that the intersection of $\Lambda_{1 2 3}$ and the flowout $\Lambda_{1 2 3}^g$ along its null directions is clean. (Note here also that a future null-geodesic starting at $K_{1 2 3}$ does not intersect $K_{1 2 3}$ again.)
%
For our purposes, $\zeta=\alpha\xi^{(0)}$ and $\zeta_j=\alpha_j\xi^{(j)}$ for some $\alpha,\alpha_j\in\mathbb{R}$.
We are particularly interested in this expression for $q=q_0$ and $y=x_0$. Notice that $q_0\in K_{123}$ and $q_0$ and $x_0$ is joined by the null-geodesic $\gamma^{(0)}$.

Now, the solution $\mathcal{U}^{(3)}$ of the initial-boundary value problem can be written as the sum of the incident wave $\mathcal{U}^{(3),\mathrm{inc}}$ and wave $\mathcal{U}^{(3),\mathrm{ref}}$ arising from reflection at $\partial M$
\[
\mathcal{U}^{(3)}=\mathcal{U}^{(3),\mathrm{inc}}+\mathcal{U}^{(3),\mathrm{ref}}.
\]
The reflected wave vanishes prior to the intersection of $\mathrm{supp}\ \mathcal{U}^{(3),\mathrm{inc}}$ with the boundary $\partial M$, and in a small neighborhood of $y$, satisfies $\square_g\mathcal{U}^{(3),\mathrm{ref}}=0$ with Neumann data $\partial_\nu\mathcal{U}^{(3),\mathrm{ref}}=-\partial_\nu\mathcal{U}^{(3),\mathrm{inc}}$. Near $y$ and in view of the null-convexity assumption on $\partial M$, the incident wave $\mathcal{U}^{(3),\mathrm{inc}}$ is a conormal distribution relative to the conormal bundle of a submanifold transversal to $\partial{M}$; therefore, so is $\mathcal{U}^{(3),\mathrm{ref}}$. Moreover, the principal symbols of the restrictions of $\mathcal{U}^{(3),\rm inc}$ and $\mathcal{U}^{(3),\rm ref}$ to $\partial M$ agree due to the Neumann boundary condition. (Indeed, following \cite{stefanov2018inverse}, we can write $\mathcal{U}^{(3),\bullet}$ in a neighborhood of $y$ in the form $\mathcal{U}^{(3),\bullet}=(2\pi)^{-3}\int e^{\mathrm{i}\phi^\bullet(x,\theta)}a^\bullet(x,\theta)\,\mathrm{d}\theta$ for $\bullet=\mathrm{inc,ref}$ and suitable symbols $a^\bullet$, where the phase functions $\phi^\bullet$ solve the eikonal equation $|\mathrm{d}\phi^\bullet|_g^2=0$ with boundary conditions $\phi^\bullet(x,\theta)=x\cdot\theta$, $x\in\partial M$, and $\partial_\nu\phi^{\mathrm{ref}}=-\partial_\nu\phi^{\mathrm{inc}}$. The Neumann boundary condition $\partial_\nu\mathcal{U}^{(3)}\vert_{\partial M}=0$ implies $(\partial_\nu\phi^{\mathrm{inc}}) a^{\mathrm{inc}}+(\partial_\nu\phi^\mathrm{ref}) a^{\mathrm{ref}}=0$, thus $a^{\mathrm{inc}}=a^{\mathrm{ref}}$ at $\partial{M}$, as claimed.)

Denote $\mathcal{R}(\mathcal{U}^{(3),\mathrm{inc}})$ to be the trace of $\mathcal{U}^{(3),\mathrm{inc}}$ on $\partial M$; the trace operator $\mathcal{R}$ an FIO of order $\frac{1}{4}$ (\cite[Chapter 5.1]{duistermaat1996fourier}) with canonical relation
\[
\Gamma_{\mathcal{R}}=\{(y_|,\eta_|,y,\eta)\in (T^*(\partial M)\times T^*M)\setminus 0; y_|=y, \eta_|=\eta\vert_{T_{y}(\partial M)}\}.
\]
For any $(y_|,\eta_|)\in T^*(\partial M)$, there exists at most one outward pointing $\eta\in L^*_yM$ such that $\eta_|=\eta\vert_{T_{y}(\partial M)}$. For such $(y_|,\eta_|,y,\eta)$, the principal symbol $\sigma^{(p)}(\mathcal{R})(y_|,\eta_|,y,\eta)$ is nonzero (cf.\ \cite[Chapter 5.1]{duistermaat1996fourier}). Using the multiplicativity of principal symbols, we then have
\begin{equation}
\label{EqDtNSymbol}
\frac{1}{2}\sigma^{(p)}\left(\partial_{\epsilon_1}\partial_{\epsilon_2}\partial_{\epsilon_3}\Lambda\left(\epsilon_1f_1+\epsilon_2f_2+\epsilon_3f_3\right)\vert_{\epsilon_1=\epsilon_2=\epsilon_3=0}\right)(y_|,\eta_|)=\sigma^{(p)}(\mathcal{R})(y_|,\eta_|,y,\eta)\sigma^{(p)}(\mathcal{U}_0^{(3),\mathrm{inc}})(y,\eta).
\end{equation}
We refer to \cite[Chapter 4]{duistermaat1996fourier} for a discussion on the compositions of FIOs.

We now show how to use this to recover $h_3$ from the principal symbol of $\mathcal{U}_0^{(3),\mathrm{inc}}$: for $j=1,2$, let $u^{(j)}$ solve the equation \eqref{maineq} with $H=H^{(j)}$ and $\partial_\nu u^{(j)}=f=\sum_{i=1}^3\epsilon_i f_i$. Decompose $\mathcal{U}^{(3),j}=\partial_{\epsilon_1}\partial_{\epsilon_2}\partial_{\epsilon_3}u^{(j)}|_{\epsilon_1=\epsilon_2=\epsilon_3=0}$ as
\[
\mathcal{U}^{(3),j}=\mathcal{U}_0^{(3),j}+\mathcal{U}_1^{(3),j},
\]
as in~\eqref{EqU3decomp}; moreover, decompose
\[
\mathcal{U}^{(3),\mathrm{inc},j}=\mathcal{U}_0^{(3),\mathrm{inc},j}+\mathcal{U}_1^{(3),\mathrm{inc},j},
\]
as in~\eqref{EqU3incDecomp}. By assumption, we have
\begin{equation}\label{Lambda_three}
\partial_{\epsilon_1}\partial_{\epsilon_2}\partial_{\epsilon_3}\Lambda^{(1)}(f)\vert_{\epsilon_1=\epsilon_2=\epsilon_3=0}=\partial_{\epsilon_1}\partial_{\epsilon_2}\partial_{\epsilon_3}\Lambda^{(2)}(f)\vert_{\epsilon_1=\epsilon_2=\epsilon_3=0};
\end{equation}
the expression~\eqref{EqDtNSymbol} shows that this implies
\begin{equation}\label{equality_three}
\sigma^{(p)}(\mathcal{U}_0^{(3),\mathrm{inc},1})(y,\eta)=\sigma^{(p)}(\mathcal{U}_0^{(3),\mathrm{inc},2})(y,\eta).
\end{equation}
By the explicit formula for $\sigma^{(p)}(\mathcal{U}_0^{(3),\mathrm{inc},j})(y,\eta)$ given by \eqref{U3}, and taking $y=x_0$, we get 
\[
h_3^{(1)}(q_0)=h_3^{(2)}(q_0).
\]
Since $q_0$ was an arbitrary point in $\mathbb{U}$, we conclude that $h_3^{(1)}=h_3^{(1)}$ in $\mathbb{U}$.

 Now we analyze
\[
\mathcal{U}_1^{(3)}:=2\sum_{\sigma\in\Sigma(3)}Q_g(h_2v_{\sigma(1)}Q_g(h_2v_{\sigma(2)}v_{\sigma(3)})).
\]
Since $h_3$ has already been recovered, we can subtract its contribution to $\mathcal{U}^{(3)}$; we can thus determine $\mathcal{U}_1^{(3)}\vert_{\partial M}$. More precisely, the fact $\partial_{\epsilon_1}\partial_{\epsilon_2}\partial_{\epsilon_3}\Lambda^{(1)}(f)\vert_{\epsilon_1=\epsilon_2=\epsilon_3=0}=\partial_{\epsilon_1}\partial_{\epsilon_2}\partial_{\epsilon_3}\Lambda^{(2)}(f)\vert_{\epsilon_1=\epsilon_2=\epsilon_3=0}$ implies
\[
(\mathcal{U}_0^{(3),1}+\mathcal{U}_1^{(3),1})\vert_{\partial M}=(\mathcal{U}_0^{(3),2}+\mathcal{U}_1^{(3),2})\vert_{\partial M}.
\]
Recall that $\mathcal{U}_0^{(3),j}:=-6Q_g(h_3^{(j)}v_1v_2v_3)$ and $h_3^{(1)}=h_3^{(2)}$ in $\mathbb{U}$; therefore, $\mathcal{U}_0^{(3),1}=\mathcal{U}_0^{(3),2}$ on $\partial M$, hence $\mathcal{U}_1^{(3),1}=\mathcal{U}_1^{(3),2}$ on $\partial M$.

 Similarly to before, we write $\mathcal{U}_1^{(3)}=\mathcal{U}_1^{(3),\mathrm{inc}}+\mathcal{U}_1^{(3),\mathrm{ref}}$, where
\[
\mathcal{U}_1^{(3),\mathrm{inc}}=2\sum_{\sigma\in\Sigma(3)}\widetilde{Q}_g(h_2v_{\sigma(1)}\widetilde{Q}_g(h_2v_{\sigma(2)}v_{\sigma(3)}))
\]
is the incident wave and $\mathcal{U}_1^{(3),\mathrm{ref}}$ is the reflected wave.  By \cite[Lemma 3.3, 3.4]{lassas2018inverse}, we have
\[
\widetilde{Q}_g(h_2v_iv_j)\in \mathcal{I}^{\mu-1,\mu}(\Lambda_{ij},\Lambda_i)+ \mathcal{I}^{\mu-1,\mu}(\Lambda_{ij},\Lambda_j).
\]
Then using \cite[Lemma 3.6 and Proposition 2.1]{lassas2018inverse}, one can obtain (cf. \cite[Proposition 3.7]{lassas2018inverse} and the discussion after it):
\begin{proposition}\label{symbol_threewaves}
For any $q\in K_{123}$ and $\zeta\in N_{q}^*K_{123}$, assume $(y,\eta)$ is joined from $(q,\zeta)$ by a null-bicharacteristic. If $h_2$ is non-vanishing on $K_{123}$, then
\[
\mathcal{U}_1^{(3),\mathrm{inc}}\in \mathcal{I}^{3\mu-\frac{3}{2},-\frac{1}{2}}(\Lambda_{123},\Lambda_{123}^g),
\]
away from $\cup_{i=1}^3\Lambda_i$, with principal symbol
\begin{align*}
\sigma^{(p)}(\mathcal{U}_1^{(3),\mathrm{inc}})(y,\eta)&=2(2\pi)^{-2}\sigma^{(p)}(\widetilde{Q}_g)(y,\eta,q,\zeta)h_2(q)^2\left(\sum_{\sigma\in\Sigma(3)}\left|\zeta_{\sigma(2)}+\zeta_{\sigma(3)}\right|_{g(q)}^{-2}\right) \\
  &\qquad \times\prod_{j=1}^3\sigma^{(p)}(v_j)(q,\zeta_j).
\end{align*}
\end{proposition}

Now we can conclude that $\sigma^{(p)}(\mathcal{U}_1^{(3),\mathrm{inc},1})(y,\eta)=\sigma^{(p)}(\mathcal{U}_1^{(3),\mathrm{inc},2})(y,\eta)$ since $\mathcal{U}_1^{(3),1}=\mathcal{U}_1^{(3),2}$ on $\partial M$; we use this for $y=x_0$ and $q=q_0$. As shown in equations~\eqref{nonvanishing30}--\eqref{nonvanishing3}, the sum $\sum_{\sigma\in\Sigma(3)}\left|\zeta_{\sigma(2)}+\zeta_{\sigma(3)}\right|_{g(q_0)}^{-2}$ appearing here is nonvanishing; therefore,
\[
(h^{(1)}_2(q_0))^2=(h^{(2)}_2(q_0))^2.
\]

\subsection{Nonlinear interactions of four waves and recovery of \texorpdfstring{$h_2$ and $h_4$}{h2 and h4}}

In this section, we use nonlinear interaction of four distorted plane waves. Thus, we take $N=4$ in \eqref{boundarysource} and consider Neumann data
\[
f=\sum_{i=1}^4\epsilon_i f_i.
\]
Take $x_1,x_2,x_3,x_4\in\widetilde{M}\setminus M$ in a neighborhood of $x_-$, where $x_-$ is as in~\eqref{Eqxminus} for some point $q_0\in\mathbb{U}$; suppose $\gamma_{x_j,\xi_j}$ joins $x_j$ to $q_0$. Take $u_i\in\mathcal{I}^\mu(\Lambda(x_i,\xi_i,s_0))$ and let $f_i=\partial_\nu u_i\vert_{\partial M}$ for $i=1,2,3,4$. One can ensure that $\Lambda_i=N^*K_i=\Lambda(x_i,\xi_i,s_0)$, $i=1,2,3,4$ satisfy Assumption \ref{assumption1} in Section \ref{notation}.

In this section, we will use the notations
\begin{gather*}
\Theta^{(1)}=\cup_{i=1}^4\Lambda_i;\quad\,\,\Theta^{(2)}=\cup_{i,j=1}^4\Lambda_{ij};\quad\,\,\Theta^{(3)}=\cup_{i,j,k=1}^4\Lambda_{ijk};\\
K^{(1)}=\cup_{i=1}^4K_i;\quad K^{(2)}=\cup_{i,j=1}^4K_{ij};\quad K^{(3)}=\cup_{i,j,k=1}^4K_{ijk}, \\
\Xi=\Theta^{(1)}\cup \Theta^{(3),g}\cup \Lambda_{q_0}.
\end{gather*}

Write
\[
\begin{split}
\mathcal{V}^{(4)}=&\partial_{\epsilon_1}\partial_{\epsilon_2}\partial_{\epsilon_3}\partial_{\epsilon_4}u\vert_{\epsilon_1=\epsilon_2=\epsilon_3=\epsilon_4=0}\\
=&-4\sum_{\sigma\in \Sigma}Q_g(h_2v_{\sigma(1)}Q_g(h_2v_{\sigma(2)}Q_g(h_2v_{\sigma(3)}v_{\sigma(4)})))\\
&-\sum_{\sigma\in \Sigma}Q_g(h_2Q_g(h_2v_{\sigma(1)}v_{\sigma(2)})Q_g(h_2v_{\sigma(3)}v_{\sigma(4)}))\\
&+2\sum_{\sigma\in \Sigma}Q_g(h_2v_{\sigma(1)}Q_g(h_3v_{\sigma(2)}v_{\sigma(3)}v_{\sigma(4)}))+3\sum_{\sigma\in \Sigma}Q_g(h_3v_{\sigma(1)}v_{\sigma(2)}Q_g(h_2v_{\sigma(3)}v_{\sigma(4)}))\\
&-24Q_g(h_4v_1v_2v_3v_4).
\end{split}
\]
Assume $\mathcal{V}^{(4)}=\mathcal{V}^{(4),\mathrm{inc}}+\mathcal{V}^{(4),\mathrm{ref}}$, where $\mathcal{V}^{(4),\mathrm{inc}}$ is the incident wave, and $\mathcal{V}^{(4),\mathrm{ref}}$ is the reflected wave.
Part of the results in \cite[Proposition 3.11, 3.12]{lassas2018inverse} can be summarized in the following proposition.
\begin{proposition}\label{symbol_fourwaves}
If $h_4(q_0)\neq 0$, we have
\[
\mathcal{V}^{(4),\mathrm{inc}}\in \mathcal{I}^{4\mu+\frac{3}{2}}(\Lambda_{q_0}^g\setminus\Xi)
\]
away from $\cup_{i=1}^3\Lambda_i$, with principal symbol
\begin{equation}\label{symbol_V4}
\sigma^{(p)}(\mathcal{V}^{(4),\mathrm{inc}})(y,\eta)=-24(2\pi)^{-3}\sigma^{(p)}(\widetilde{Q}_g)(y,\eta,q_0,\zeta)h_4(q_0)\prod_{j=1}^4\sigma^{(p)}(v_j)(q_0,\zeta_j),
\end{equation}
for $(y,\eta)\in\Lambda_{q_0}^g\setminus\Xi$.
Here $(y,\eta)$ is joined with $(q_0,\zeta)$ by a null-bicharacteristic of $\square_g$, and $\zeta\in L_{q_0}^{*,+}M$ has the unique decomposition $\zeta=\sum_{i=4}^4\zeta_i$ with $\zeta_i\in N_{q_0}^*K_i$.
\end{proposition}

Assume $h_3^{(1)},h_4^{(2)}\neq 0$ at $q_0$. Denote $\mathcal{K}^{(3)}=\pi(\Theta^{(3),g})\subset M$. By taking $s_0\rightarrow 0$, the set $K^{(1)}\cup\mathcal{K}^{(3)}$ tends to a set of Hausdorff dimension 2 (cf. \cite[Section 4]{kurylev2018inverse}). Thus we can choose $s_0$ small enough such that there exists $\zeta\in\Lambda_{q_0}\setminus(\Theta^{(1)}\cup \Theta^{(3)})$ such that $y\in(0,T)\times \partial N$. But then
\[
\partial_{\epsilon_1}\partial_{\epsilon_2}\partial_{\epsilon_3}\partial_{\epsilon_4}\Lambda^{(1)}(f)\vert_{\epsilon_1=\epsilon_2=\epsilon_3=\epsilon_4=0}=\partial_{\epsilon_1}\partial_{\epsilon_2}\partial_{\epsilon_3}\partial_{\epsilon_4}\Lambda^{(2)}(f)\vert_{\epsilon_1=\epsilon_2=\epsilon_3=\epsilon_4=0}
\]
implies
 \[
 \sigma^{(p)}(\mathcal{V}^{(4),\mathrm{inc},1})(y,\eta)=\sigma^{(p)}(\mathcal{V}^{(4),\mathrm{inc},2})(y,\eta).
 \]
 By the explicit expression for $\sigma^{(p)}(\mathcal{V}^{(4),\mathrm{inc},j})(y,\eta)$ given in \eqref{symbol_V4}, we obtain
\[
h_4^{(1)}(q_0)=h_4^{(2)}(q_0).
\]

With $h_4$ thus recovered in $\mathbb{U}$, we can determine
\[
\mathcal{V}^{(4)}_1=\mathcal{V}^{(4)}+24Q_g(h_4v_1v_2v_3v_4).
\]
at the boundary $ (0,T)\times\partial N$. Here we use the fact that, by the finite speed of propagation, $Q_g(h_4v_1v_2v_3v_4)\vert_{(0,T)\times\partial N}$ depends only on the value of $h_4v_1v_2v_3v_4$ in $J^-((0,T)\times\partial N)$ and $v_j$ vanishes on $M\setminus J^+((0,T)\times\partial N)$.
Similar as the previous section, we can write $\mathcal{V}^{(4)}_1=\mathcal{V}^{(4),\mathrm{inc}}_1+\mathcal{V}^{(4),\mathrm{ref}}_1$, which is the sum of the incident wave and reflected wave. The microlocal property of $\mathcal{V}^{(4),\mathrm{inc}}_1$ is analyzed carefully in the proofs of \cite[Proposition 3.11, 3.12]{lassas2018inverse}. We summarize the results that we need in the following proposition.
\begin{proposition}\label{symbol_fourwaves}
Assume $(y,\eta)\in \Lambda_{q_0}^g\setminus \Xi$ is joined from $(q_0,\zeta)\in\Lambda_{q_0}$ by a null-bicharacteristic. 
\begin{enumerate}
\item If $h_3(q_0)\neq 0$, we have $\mathcal{V}^{(4),\mathrm{inc}}_1\in \mathcal{I}^{4\mu-\frac{1}{2}}(\Lambda_{q_0}^g\setminus\Xi)$ with principal symbol
\[
\begin{split}
\sigma^{(p)}(\mathcal{V}^{(4),\mathrm{inc}}_1)(y,\eta)=(2\pi)^{-3}h_2(q_0)h_3(q_0)\mathcal{G}_2(\zeta)\sigma^{(p)}(Q_g)(y,\eta,q_0,\zeta)\prod_{j=1}^4\sigma^{(p)}(v_j)(q_0,\zeta_j),
\end{split}
\]
where
\[
\mathcal{G}_2(\zeta)=\sum_{\sigma\in\Sigma(4)}\left(\frac{3}{|\zeta_{\sigma(1)}+\zeta_{\sigma(2)}|^2_{g(q_0)}}+\frac{2}{|\zeta_{\sigma(2)}+\zeta_{\sigma(3)}+\zeta_{\sigma(4)}|^2_{g(q_0)}}\right).
\]
\item If $h_3= 0$ in a neighborhood of $q_0$, we have $\mathcal{V}^{(4),\mathrm{inc}}_1\in \mathcal{I}^{4\mu-\frac{5}{2}}(\Lambda_{q_0}^g\setminus\Xi)$ with principal symbol
\[
\begin{split}
\sigma^{(p)}(\mathcal{V}^{(4),\mathrm{inc}}_1)(y,\eta)=(2\pi)^{-3}h_2(q_0)^3\mathcal{G}_3(\zeta)\sigma^{(p)}(Q_g)(y,\eta,q_0,\zeta)\prod_{j=1}^4\sigma^{(p)}(v_j)(q_0,\zeta_j),
\end{split}
\]
where
\[
\mathcal{G}_3(\zeta)=\sum_{\sigma\in\Sigma(4)}\left(\frac{4}{|\zeta_{\sigma(2)}+\zeta_{\sigma(3)}+\zeta_{\sigma(4)}|^2_{g(q_0)}}+\frac{1}{|\zeta_{\sigma(1)}+\zeta_{\sigma(2)}|^2_{g(q_0)}}\right)\frac{1}{|\zeta_{\sigma(3)}+\zeta_{\sigma(4)}|^2_{g(q_0)}}.
\]
\end{enumerate}
\end{proposition}

Now $\Lambda^{(1)}=\Lambda^{(2)}$ implies 
\[
\sigma^{(p)}(\mathcal{V}^{(4),\mathrm{inc},1}_1)(y,\eta)=\sigma^{(p)}(\mathcal{V}^{(4),\mathrm{inc},2}_1)(y,\eta).
\]
Using Proposition \ref{symbol_fourwaves}, and the (generic) nonvanishing of $\mathcal{G}_2$ and $\mathcal{G}_3$ (\cite[Proposition 3.12]{lassas2018inverse}), we now have
\[
h_2^{(1)}(q_0)h_3^{(1)}(q_0)=h_2^{(2)}(q_0)h_3^{(2)}(q_0)
\]
if $h_3^{(j)}(q_0)\neq 0$ or
\begin{equation}
\label{h2cubed}
h_2^{(1)}(q_0)^3=h_2^{(2)}(q_0)^3.
\end{equation}
if $h_3^{(j)}$ vanishes near $q_0$. For either case, we can obtain
\[
h_2^{(1)}(q_0)=h_2^{(2)}(q_0),
\]
invoking the facts $h_2^{(1)}(q_0)^2=h_2^{(2)}(q_0)^2$ and $h_3^{(1)}(q_0)=h_3^{(2)}(q_0)$. If $h_3^{(j)}$ vanishes at $q_0$ but not nearby, then we are in case~\eqref{h2cubed} at a sequence of points tending to $q_0$, hence obtaining the equality $h_3^{(1)}(q_0)=h_3^{(2)}(q_0)=0$ by continuity.

\subsection{Recovery of \texorpdfstring{$h_k$, $k\geq 5$}{hk for k at least 5}}
Finally, we recover $h_k$ for $k=5,6,\ldots,$ using the interaction of three waves. The coefficients $h_2,h_3,h_4$ have already been determined above. Inductively, assume that all $h_k$, $k\leq N-1$ $(N\geq 5)$, have already be recovered; we proceed to recover $h_N$. Denote
\[
\mathcal{U}^{(N)}=\partial^{N-2}_{\epsilon_1}\partial_{\epsilon_2}\partial_{\epsilon_3}u\vert_{\epsilon_1=\epsilon_2=\epsilon_3=0},
\]
where $u$ is the solution to \eqref{maineq} with $f=\sum_{i=1}^3\epsilon_if_i$. We observe that
\[
\mathcal{U}^{(N)}=-N!Q_g(h_Nv_1^{N-2}v_2v_3)+R_N(v_1,v_2,v_3;h_2,\ldots,h_{N-1}),
\]
where $R_N(v_1,v_2,v_3;h_2,\ldots,h_{N-1})$ depends on $v_1,v_2,v_3$ and $h_2,\ldots,h_{N-1}$ only. We note here that the singularities in $R_N$ are very complicated. The Sobolev regularity of $R_N$ was analyzed in \cite[Section 5]{lassas2018inverse} on boundaryless Lorentzian manifolds. We avoid the complication by using the following inductive procedure.

 Now, $h_2,\ldots,h_{N-1}$ have already been recovered in $\mathbb{U}$; moreoever, $v_1,v_2,v_3$ (which vanish on $M\setminus J^+((0,T)\times\partial N)$) are known; hence, $R_N$ is known on $(0,T)\times\partial N$ by finite speed of propagation.
 Thus we can recover
\[
\mathcal{U}^{(N)}_0=-N!Q_g(h_Nv_1^{N-2}v_2v_3)
\]
on the boundary $(0,T)\times\partial N$ from $\Lambda$. Assume $\mathcal{U}^{(N)}_0=\mathcal{U}^{(N),\mathrm{inc}}_0+\mathcal{U}^{(N),\mathrm{ref}}_0$, where
\[
\mathcal{U}_0^{(N),\mathrm{inc}}=-N!\widetilde{Q}_g(h_Nv_1^{N-2}v_2v_3).
\]
By \cite[Lemma 5.1]{lassas2018inverse}, we have $v_1^{N-2}\in\mathcal{I}^{\mu+(N-3)(\mu+\frac{3}{2})}(\Lambda_1)$, with
\[
\sigma^{(p)}(v_1^{N-2})=(2\pi)^{-\frac{N-3}{2}}\underbrace{\sigma^{(p)}(v_1)*\sigma^{(p)}(v_1)*\cdots *\sigma^{(p)}(v_1)}_{N-2\text{ factors},\ N-3\text{ convolutions}}=:(2\pi)^{-\frac{N-3}{2}}A_1^{(N-2)}.
\]
By the proof of \cite[Proposition 5.6]{lassas2018inverse}, $A_1^{(N-2)}$ is non-vanishing at $(q_0,\zeta_1)$. By \cite[Lemma 3.3]{lassas2018inverse}, $v_2v_3\in\mathcal{I}^{\mu,\mu+1}(\Lambda_{23},\Lambda_2)+\mathcal{I}^{\mu,\mu+1}(\Lambda_{23},\Lambda_3)$, and then by \cite[Lemma 3.6]{lassas2018inverse}
\[
v_1^{N-2}v_2v_3\in \mathcal{I}^{3\mu+(N-3)(\mu+\frac{3}{2})}(\Lambda_{123})\quad \text{away from $\cup_{i=1}^3\Lambda_i$}.
\]
By \cite[Proposition 2.1]{lassas2018inverse}, {we have
\begin{proposition}\label{symbol_threewaves1}
If $h_N$ is non-vanishing on $K_{123}$, we have
\[
\mathcal{U}_0^{(N),\mathrm{inc}}\in\mathcal{I}^{3\mu+(N-3)(\mu+\frac{3}{2})+\frac{1}{2},-\frac{1}{2}}(\Lambda_{123},\Lambda_{123}^g),
\]
away from $\cup_{i=1}^3\Lambda_i$, with principal symbol
\begin{equation}\label{Uk}
\begin{split}
&\sigma^{(p)}(\mathcal{U}_0^{(N),\mathrm{inc}})(y,\eta)\\
&\qquad = -N!(2\pi)^{-2-\frac{N-3}{2}}\sigma^{(p)}(\widetilde{Q}_g)(y,\eta,q_0,\zeta)h_N(q_0)A_1^{(N-2)}(q_0,\zeta_1)\prod_{j=2}^3\sigma^{(p)}(v_j)(q_0,\zeta_j).
\end{split}
\end{equation}
\end{proposition}
}

As around \eqref{Lambda_three} and \eqref{equality_three} (and using the same notation), the equality
\[
\partial^{N-2}_{\epsilon_1}\partial_{\epsilon_2}\partial_{\epsilon_3}\Lambda^{(1)}(f)\vert_{\epsilon_1=\epsilon_2=\epsilon_3=0}=\partial^{N-2}_{\epsilon_1}\partial_{\epsilon_2}\partial_{\epsilon_3}\Lambda^{(2)}(f)\vert_{\epsilon_1=\epsilon_2=\epsilon_3=0}
\]
thus implies
 \[
 \sigma^{(p)}(\mathcal{U}_0^{(N),\mathrm{inc},1})(y,\eta)=\sigma^{(p)}(\mathcal{U}_0^{(N),\mathrm{inc},2})(y,\eta).
 \]
  By the explicit formula for $\sigma^{(p)}(\mathcal{U}_0^{(N),\mathrm{inc},j})(y,\eta)$ given by \eqref{Uk}, we get 
  \[
  h_N^{(1)}(q_0)=h_N^{(2)}(q_0).
  \]
This completes the proof of Theorem \ref{maintheorem}.
\section{Recovery using Gaussian beams}\label{gaussianbeam}

In this section, we give an alternative approach to recover $H$, assuming $h_2$ is \textit{a priori} known, using Gaussian beam solutions to the linear wave equation. Such approach for nonlinear wave equations have been undertaken in \cite{kurylev2013determination,feizmohammadi2019recovery,oksanen2020inverse}. We note here that Gaussian beams have also been used for various inverse problems \cite{katchalov1998multidimensional,bao2014sensitivity, belishev1996boundary, dos2016calderon,feizmohammadi2019timedependent,feizmohammadi2019recovery,feizmohammadi2019inverse}.

We still use higher order linearization of the Neumann-to-Dirichlet map $\Lambda$, but will obtain an integral identity and use it to recover the parameters. Gaussian beams will be used in the integral identity. A similar technique was applied to a nonlinear elastic wave equation in \cite{uhlmann2019inverse}. Higher order linearizations of the Dirichlet-to-Neumann map and the resulting integral identities for semilinear and quasilinear elliptic equations have been used in \cite{sun1997inverse,kang2002identification,assylbekov2017direct,carstea2019reconstruction,lassas2019inverse,lassas2019partial,feizmohammadi2019inverse,krupchyk2020remark,krupchyk2019partial}.

Let $v_j$, $j=1,2,\ldots$, solve
\begin{alignat}{2}
\square_g v_j&=0&\quad&\text{ in }(0,T)\times N,\nonumber\\
\label{eq_vj}
\partial_\nu v_j&=f_j&\quad&\text{ on }(0,T)\times \partial N,\\
v_j=\partial_tv_j&=0&\quad&\text{ on }\{t=0\}.\nonumber
\end{alignat}
Let $v_0$ be the solution to the backward wave equation
\begin{alignat}{2}
\square v_0&=0&\quad&\text{ in }(0,T)\times  N,\nonumber\\
\label{eq_v0}
\partial_\nu v_0&=f_0&\quad&\text{ on }(0,T)\times \partial N,\\
v_0=\partial_tv_0&=0&\quad&\text{ on }\{t=T\}.\nonumber
\end{alignat}

First let us recover $h_3$. Take $f=\epsilon_1f_1+\epsilon_2f_2+\epsilon_3f_3$, and let $u$ solve \eqref{maineq}. Denote $\mathcal{U}^{(123)}=\frac{\partial^3}{\partial\epsilon_1\partial\epsilon_2\partial\epsilon_3}u\vert_{\epsilon_1=\epsilon_2=\epsilon_3=0}$, $\mathcal{U}^{(ij)}=\frac{\partial^2}{\partial\epsilon_i\partial\epsilon_j}u\vert_{\epsilon_i=\epsilon_j=0}$. Notice $\frac{\partial}{\partial\epsilon_i}u\vert_{\epsilon_i=0}=v_i$ and $\mathcal{U}^{(ij)}$ solves
\begin{alignat*}{2}
\square\mathcal{U}^{(ij)}+h_2(x)v_iv_j&=0&\quad&\text{ in }(0,T)\times  N\\
\partial_\nu\mathcal{U}^{(ij)}&=0&\quad&\text{ on }(0,T)\times \partial N,\\
\mathcal{U}^{(ij)}=\partial_t\mathcal{U}^{(ij)}&=0&\quad&\text{ on }\{t=0\}.
\end{alignat*}
 Applying $\frac{\partial^3}{\partial\epsilon_1\partial\epsilon_2\partial\epsilon_3}$ to \eqref{maineq} evaluated at at $\epsilon_1=\epsilon_2=\epsilon_3=0$, we get
\[
\square \mathcal{U}^{(123)}+h_2(x)\sum_{\sigma\in\Sigma(3)}\mathcal{U}^{(\sigma(1)\sigma(2))}v_{\sigma(3)}+6h_3(x)v_1v_2v_3=0.
\]
Integration by parts gives
\begin{align}
&\int_{\partial M}\frac{\partial^3}{\partial\epsilon_1\partial\epsilon_2\partial\epsilon_3}\Big\vert_{\epsilon_1=\epsilon_2=\epsilon_3=0}\Lambda(\epsilon_1f_1+\epsilon_2f_2+\epsilon_3f_3)f_0\,\mathrm{d}V_g\nonumber\\
\label{IBP}
&\qquad=\int_Mh_3v_1v_2v_3v_0\,\mathrm{d}V_g+\int_Mh_2(x)\sum_{\sigma\in\Sigma(3)}\mathcal{U}^{(\sigma(1)\sigma(2))}v_{\sigma(3)} v_0\,\mathrm{d}V_g.
\end{align}
we note here that by finite speed of propagation for solutions of the wave equation, the functions $v_i,v_j$ and thus also $\mathcal U^{(i j)}$ vanish in $M\setminus J^+((0,T)\times\partial N)$, $i,j=1,2,3$, and likewise $v_0$ vanishes in $M\setminus J^-((0,T)\times\partial N)$; therefore, our knowledge of $h_2$ in $\mathbb{U}$ is sufficient to compute the second summand in~\eqref{IBP}. Therefore, we can recover
\begin{equation}\label{integral_h3}
\int_Mh_3v_1v_2v_3v_0\,\mathrm{d}V_g.
\end{equation}
We will use special solutions $v_1,v_2,v_3,v_0$ in the above identity and thereby recover the coefficient $h_3$. Concretely, we shall use Gaussian beam solutions for the wave equation $\square_g v=0$ on $\widetilde M$ of the form
\[
v(x)=e^{\mathrm{i}\rho\varphi(x)}\mathfrak{a}_\rho(x)+R_\rho(x),
\]
with a large parameter $\rho$. The phase function $\varphi$ is complex-valued. The principal term $e^{\mathrm{i}\rho\varphi(x)}\mathfrak{a}(x)$ is concentrated near a null geodesic $\gamma$ in the manifold $\mathbb{R}\times N$. 
The remainder term $R_\rho$ will vanish rapidly as $\rho\rightarrow+\infty$.

\textbf{Fermi coordinates on $\widetilde M$.}  Assume $\gamma$ passes through a point $p\in M$ and joins two points $\gamma(\tau_-)$ and $\gamma(\tau_+)$ on the boundary $\mathbb{R}\times\partial N$. We will use the Fermi coordinates $\Phi$ on $\widetilde M$ in a neighborhood of $\gamma([\tau_-,\tau_+])$, denoted by $(z^0:=\tau,z^1,z^2,z^3)$, such that $\Phi(\gamma(\tau))=(\tau,0)$ (cf.\ \cite[Lemma 1]{feizmohammadi2019recovery}). 

\textbf{Construction of Gaussian beams.} We will construct asymptotic solutions of the form $u_\rho=\mathfrak{a}_\rho e^{\mathrm{i}\rho\varphi}$ on $\widetilde M$ with 
\[
\varphi=\sum_{k=0}^N\varphi_k(\tau,z'),\quad \mathfrak{a}_\rho(\tau,z')=\chi\left(\frac{|z'|}{\delta}\right)\sum_{k=0}^N\rho^{-k}a_k(\tau,z'),\quad
a_k(\tau,z')=\sum_{j=0}^Na_{k,j}(\tau,z')
\]
in a neighborhood of $\gamma$,
\begin{equation}\label{neighbor_V}
\mathcal{V}=\bigl\{(\tau,z')\in\widetilde{M} : \tau\in \bigl[\tau_--\tfrac{\epsilon}{\sqrt{2}},\tau_++\tfrac{\epsilon}{\sqrt{2}}\bigr],\,|z'|<\delta\bigr\}.
\end{equation}
Here for each $j$, $\varphi_j$ and $a_{k,j}$ are a complex valued homogeneous polynomials of degree $j$ with respect to the variables $z^i$, $i=1,2,3$, and $\delta>0$ is a small parameter. The smooth function $\chi:\mathbb{R}\rightarrow [0,+\infty)$ satisfies $\chi(t)=1$ for $|t|\leq\frac{1}{4}$ and $\chi(t)=0$ for $|t|\geq \frac{1}{2}$. 

We have
\begin{equation}
\begin{split}
\square_g(\mathfrak{a}_\rho e^{\mathrm{i}\rho\varphi})&=e^{\mathrm{i}\rho\varphi}(\rho^2(\mathcal{S}\varphi)\mathfrak{a}_\rho-\mathrm{i}\rho\mathcal{T}\mathfrak{a}_\rho+\square_g\mathfrak{a}_\rho), \\
  &\quad \mathcal{S}\varphi=\langle \mathrm{d}\varphi,\mathrm{d}\varphi\rangle_g, \\
  &\quad \mathcal{T}a=2\langle \mathrm{d}\varphi,\mathrm{d}a\rangle_g-(\square_g\varphi)a.
\end{split}
\end{equation}
We need to construct $\varphi$ and $\mathfrak{a}_\rho$ such that
\begin{equation}\label{S_cond}
\frac{\partial^\Theta}{\partial z^\Theta}(\mathcal{S}\varphi)(\tau,0)=0,\quad
\frac{\partial^\Theta}{\partial z^\Theta}(\mathcal{T}a_0)(\tau,0)=0, \quad
\frac{\partial^\Theta}{\partial z^\Theta}(-\mathrm{i}\mathcal{T}a_k+\square_g a_{k-1})(\tau,0)=0
\end{equation}
for $\Theta=(0,\Theta_1,\Theta_2,\Theta_3)$ with $|\Theta|\leq N$. For more details we refer to \cite{feizmohammadi2019recovery}. Following \cite{feizmohammadi2019timedependent}, we take
\[
\varphi_0=0,\quad \varphi_1=z^1,\quad
\varphi_2(\tau,z)=\sum_{1\leq i,j\leq 3}H_{ij}(\tau)z^iz^j.
\]
Here $H$ is a symmetric matrix with $\Im H(\tau)>0$; the matrix $H$ satisfies a Riccati ODE,
\begin{equation}\label{Ricatti}
\frac{\mathrm{d}}{\mathrm{d}\tau}H+HCH+D=0,\quad\tau\in \bigl(\tau_--\tfrac{\epsilon}{2},\tau_++\tfrac{\epsilon}{2}\bigr),\quad H(0)=H_0,\text{ with }\Im H_0>0,
\end{equation}
where $C$, $D$ are matrices with $C_{11}=0$, $C_{ii}=2$, $i=2,3$, $C_{ij}=0$, $i\neq j$ and $D_{ij}=\frac{1}{4}(\partial_{ij}^2g^{11})$.
 
 \begin{lemma}[\text{\cite[Lemma 3.2]{feizmohammadi2019timedependent}}]
 The Ricatti equation \eqref{Ricatti} has a unique solution. Moreover the solution $H$ is symmetric and $\Im (H(\tau))>0$ for all $\tau\in (\tau_--\frac{\delta}{2},\tau_++\frac{\delta}{2})$. For solving the above Ricatti equation, one has $H(\tau)=Z(\tau)Y(\tau)^{-1}$, where $Y(\tau)$ and $Z(\tau)$ solve the ODEs
 \[
 \begin{split}
& \frac{\mathrm{d}}{\mathrm{d}\tau}Y(\tau)=CZ(\tau),\quad Y(0)=Y_0,\\
&\frac{\mathrm{d}}{\mathrm{d}\tau}Z(\tau)=-D(\tau)Y(\tau),\quad Z(0)=Y_1=H_0Y_0.
 \end{split}
 \]
 In addition, $Y(\tau)$ is nondegenerate.
 \end{lemma}

 \begin{lemma}[\text{\cite[Lemma 3.3]{feizmohammadi2019timedependent}}]\label{lemma_H0}
 The following identity holds:
\[
 \det(\Im(H(\tau))|\det(Y(\tau))|^2=c_0
\]
 with $c_0$ independent of $\tau$.
 \end{lemma}
  We see that the matrix $Y(\tau)$ satisfies
  \begin{equation}\label{eq_Y}
  \frac{\mathrm{d}^2}{\mathrm{d}\tau^2}Y+CD Y=0,\quad Y(0)=Y_0,\quad \frac{\mathrm{d}}{\mathrm{d}\tau}Y(0)=CY_1.
  \end{equation}

As in \cite{feizmohammadi2019recovery}, we have the following estimate by the construction of $u_\rho$ (cf.\ \eqref{S_cond})
\begin{equation}\label{est_remainder}
\|\square_g u_\rho\|_{H^k(M)}\leq C\rho^{-K},\qquad K=\frac{N+1-k}{2}-1.
\end{equation}

Consider a point $p\in\mathbb{U}$, let $x_j$, $j=0,1,2,3$ be the points on $(0,T)\times N$ chosen in Section \ref{threewaves}, and $\gamma^{(j)}$ the null-geodesics passing through $x_j$ and $q_0$. The null-geodesic $\gamma^{(j)}$ could not be self-intersecting by the global hyperbolicity of $M$. Also $\xi^{(j)}\in L^{*,+}_{q_0}M$ is the cotangent vector to $\gamma^{(j)}$ at $q_0$. By the discussions in Section \ref{threewaves}, there exits constant $\kappa_j$, $j=0,1,2,3$ such that
\begin{equation}\label{fourvectors}
\kappa_0 \xi^{(0)}+\kappa_1 \xi^{(1)}+\kappa_2 \xi^{(2)}+\kappa_3 \xi^{(3)}=0.
\end{equation}
We construct Gaussian beams $u^{(j)}_\rho$, $j=0,1,2,3$ as above of the form
\[
u^{(j)}_\rho=e^{\mathrm{i}\kappa_j\rho\varphi^{(j)}}\mathfrak{a}^{(j)}_{\kappa_j\rho},
\]
which is compactly supported in the neighborhood $\mathcal{V}$ of the null-geodeisc $\gamma^{(j)}$ (cf.\ \eqref{neighbor_V}). The parameter $\delta$ can be taken small enough such that $u_\rho^{(j)}=0$ near $\{t=0\}$ for $j=1,2,3$ and $u_\rho^{(0)}=0$ near $\{t=T\}$.

For $j=1,2,3$, we can construct a solution $v_j$ for the initial boundary value problem \eqref{eq_vj} of the form $v_j=u_\rho^{(j)}+R_\rho^{(j)}$, where the remainder term $R_\rho^{(1)}$ is a solution of
\begin{alignat*}{2}
\square_g R_\rho^{(j)}&=-\square_g u_\rho^{(1)}&\quad&\text{ on }\partial N\times (0,T),\\
\partial_\nu R_\rho^{(j)}&=0&\quad&\text{ on } \partial N\times (0,T),\\
R_\rho^{(j)}=\partial_tR_\rho^{(j)}&=0&\quad&\text{ on }\{t=0\}.
\end{alignat*}
We note here that $v_j=u_\rho^{(j)}+R_\rho^{(j)}$ is the solution to \eqref{eq_vj} with boundary value $f_j=\partial_\nu u_\rho^{(j)}\vert_{\partial M}$.
Invoking \eqref{est_remainder}, the solution $R_\rho^{(j)}$ satisfies the estimate (cf. \cite[Theorem 3.1]{dafermos1985energy}, \cite[Proposition 2.2]{oksanen2020inverse})
\[
\|R_\rho^{(j)}\|_{H^{k+1}(M)}\leq C\rho^{-K}.
\]
Using Sobolev embedding, we can choose $N$ large enough such that
\begin{equation}\label{est_R}
\|R_\rho^{(j)}\|_{C(M)}\leq C\rho^{-\frac{n+1}{2}-2}.
\end{equation}

Similarly, we can construct a solution to \eqref{eq_v0} of the form $v_0=u^{(0)}_\rho+R^{(0)}_\rho$. We only need to take the remainder term $R_\rho^{(0)}$ to be the solution to the initial value problem
\[
\begin{split}
&\square_g R_\rho^{(0)}=-\square_g u_\rho^{(0)}\\
&\partial_\nu R_\rho^{(0)}=0\text{ on }\partial N\times (0,T),\\
&R_\rho^{(0)}=\partial_tR_\rho^{(0)}=0\text{ on }\{t=0\}.
\end{split}
\]
Now $v_0$ is the solution to \eqref{eq_v0} with $g=\partial_\nu u_\rho^{(0)}\vert_{\partial M}$.


Then by the estimate \eqref{est_R}, the Neumann-to-Dirichlet map determines
\begin{equation}\label{integral}
\begin{split}
\mathcal{I}=&\rho^{\frac{n+1}{2}}\int_Mh_3v_1v_2v_3v_0\,\mathrm{d}V_g\\
=&\rho^{\frac{n+1}{2}}\int_Mh_3e^{\mathrm{i}\rho(\kappa_0\varphi^{(0)}+\kappa_1\varphi^{(1)}+\kappa_2\varphi^{(2)}+\kappa_3\varphi^{(3)})}\mathfrak{a}^{(0)}_{\kappa_0\rho}\mathfrak{a}^{(1)}_{\kappa_1\rho}\mathfrak{a}^{(2)}_{\kappa_3\rho}\mathfrak{a}^{(3)}_{\kappa_3\rho}\,\mathrm{d}V_g+\mathcal{O}(\rho^{-1}).
\end{split}
\end{equation}

 \begin{lemma}[\text{\cite[Lemma 5]{feizmohammadi2019recovery}}]
 The function
 \[
 S:=\kappa_0\varphi^{(0)}+\kappa_1\varphi^{(1)}+\kappa_2\varphi^{(2)}+\kappa_3\varphi^{(3)}
 \]
 is well-defined in a neighborhood of $q_0$ and
 \begin{enumerate}
 \item $S(q_0)=0$;
 \item $\nabla S(q_0)=0$;
 \item $\Im S(q)\geq cd(q,q_0)^2$ for $q$ in a neighborhood of $q_0$, where $c>0$ is a constant.
 \end{enumerate}
 \end{lemma}
The four null-geodesics $\gamma^{(j)}$, $j=0,1,2,3$ intersect only at the point $q_0$, invoking the condition that cut points do not exist. Therefore the product $\mathfrak{a}^{(0)}_{\kappa_0\rho}\mathfrak{a}^{(1)}_{\kappa_1\rho}\mathfrak{a}^{(2)}_{\kappa_3\rho}\mathfrak{a}^{(3)}_{\kappa_3\rho}$ is supported in a neighborhood of $q_0$. By the above lemma, and applying stationary phase (cf., for example, \cite[Theorem 7.7.5]{hormander2015analysis}) to \eqref{integral}, we have
 \[
 c\mathcal{I}=h_3(q_0)a^{(0)}_0(q_0)a^{(1)}_0(q_0)a^{(2)}_0(q_0)a^{(3)}_0(q_0)+\mathcal{O}(\rho^{-1}),
 \]
 for some explicit constant $c\neq 0$. Hence the Neumann-to-Dirichlet map $\Lambda$ determines $h_3(q_0)$.

Next we recover the higher order coefficients $h_k$, $k=4,5,\ldots$. Recursively, assume we have already recovered $h_3,\ldots, h_{N-1}$, $N\geq 4$, in $\mathbb{U}$. To recover $h_N$, take $f=\sum_{k=1}^N\epsilon_kf_k$ and apply $\frac{\partial^N}{\partial\epsilon_1\cdots\partial\epsilon_N}$ to \eqref{maineq} evaluated at at $\epsilon_1=\cdots=\epsilon_N=0$, we get the equation for $\mathcal{U}^{(12\cdots N)}=\frac{\partial^N}{\partial\epsilon_1\cdots\partial\epsilon_N}u$
\[
\begin{split}
\square\mathcal{U}^{(12\cdots N)}+R_N(v_1,\ldots, v_N;h_1,\ldots, h_{N-1})+N!h_N\prod_{k=1}^Nv_k&=0\text{ in }N\times (0,T),\\
\partial_\nu\mathcal{U}^{(12\cdots N)}&=0\text{ on }\partial N\times (0,T).
\end{split}
\]
By the recursive assumption, $R_N(v_1,\ldots, v_N,h_1,\ldots, h_{N-1})$ is already known. By integration by parts, we have
\[
\begin{split}
&\int_{\partial M}\frac{\partial^N}{\partial\epsilon_1\cdots\partial\epsilon_N}\Big\vert_{\epsilon_1=\cdots=\epsilon_N=0}\Lambda\left(\sum_{k=1}^N\epsilon_kf_k\right)g\,\mathrm{d}S_g\\
=&\int_MN!h_Nv_1\cdots v_Nv_0\,\mathrm{d}V_g+\int_M R_N(v_1,\ldots, v_N;h_1,\ldots, h_{N-1}) v_0\,\mathrm{d}V_g.
\end{split}
\]
Thus, we can recover
\begin{equation}\label{integral_hN}
\int_Mh_Nv_0v_1\cdots v_N\,\mathrm{d}V_g.
\end{equation}
Take
\[
\begin{split}
u^{(0)}_\rho&=e^{\mathrm{i}\kappa_0\rho\varphi^{(0)}}\mathfrak{a}^{(0)}_{\kappa_0\rho},\\
u^{(j)}_\rho&=e^{\mathrm{i}\kappa_j\rho\varphi^{(j)}}\mathfrak{a}^{(j)}_{\kappa_j\rho},\quad j=1,2,\\
u^{(j)}_\rho&=e^{\mathrm{i}\frac{\kappa_3}{N-2}\rho\varphi^{(3)}}\mathfrak{a}^{(3)}_{\frac{\kappa_3}{N-2}\rho},\quad j=3,\ldots, N.
\end{split}
\]
Take $f_j=\partial_\nu v^{(j)}\vert_{\partial M}$, $j=1,\ldots,N$, $g=\partial_\nu v^{(0)}_\rho\vert_{\partial M}$ this time. Then we can recover
\[
\rho^{\frac{n+1}{2}}\int_Mh_Ne^{\mathrm{i}\rho S_0}\mathfrak{a}^{(0)}_{\kappa_0\rho}\mathfrak{a}^{(1)}_{\kappa_1\rho}\mathfrak{a}^{(2)}_{\kappa_2\rho}(\mathfrak{a}^{(3)}_{\frac{\kappa_3}{N-2}\rho})^{N-2}\,\mathrm{d}V_g+\mathcal{O}(\rho^{-1}).
\]
Again applying stationary phase, we can recover $h_N(q_0)$.

\section{Discussion}\label{discussion}

We can see that $h_2$ is more difficult to recover than $h_k$, $k=3,4,\ldots$. Indeed, we need to exploit the interaction of four waves (associated with four future light-like vectors) in Section \ref{distorted}; three light-like vectors are not sufficient. (And certainly not two: as pointed out in \cite{lassas2018inverse}, the interaction of two conormal waves does not produce new propagating singularities.)

The use of Gaussian beams avoids some involved microlocal analysis and simplifies the proof substantially. 
 In our problem, we are however unable to recover $h_2$ using Gaussian beams. 
  Despite their difference, the two approaches recover $h_k$ for $k\geq 3$ in a very similar way. They both choose solutions $v_1,\ldots, v_k$ such that $v_1v_2\cdots v_k$ is supported in a neighborhood of a single point $q_0\in\mathbb{U}$ at which one wishes to determine $h_k(q_0)$.

Distorted plane waves and Gaussian beams can be constructed even when conjugate points exist. In this paper, we assume that conjugate points do not exist for the sake of simplicity of exposition. Since we prove that local recovery is possible, a layer stripping strategy as used in \cite{kurylev2018inverse} can be applied if there are conjugate points. We note that the article \cite{feizmohammadi2019recovery} determines a zeroth order potential using nonlinear interactions and does allow for the presence of conjugate points.

\subsection*{Acknowledgements}
The authors are very grateful to two careful referees for their detailed and helpful suggestions and corrections. GU was partially supported by NSF, a Walker Professorship at UW and a Si-Yuan Professorship at IAS, HKUST. PH, GU as a senior Clay Scholar, and JZ acknowledge the great hospitality of MSRI, where part of this work was carried out during their visits. Part of this research was conducted during the period PH served as a Clay Research Fellow.

\end{document}